\documentclass[12pt]{article}
\usepackage{amsmath}
\usepackage{amssymb}

\newcommand{\R}{{\mathbb R}}
\newcommand{\N}{{\mathbb N}}

\newcommand{\fn}{\!:\!}
\newcommand{\lsum}{\sum\limits}
\newcommand{\llim}{\lim\limits}
\newcommand{\lint}{\int\limits}
\newcommand{\qed}{\mbox{$\quad\blacksquare$}}
\newtheorem{theorem}{Theorem}
\newtheorem{corollary}[theorem]{Corollary}
\begin{document}
\begin{center}
{\large\bf Necessary and sufficient conditions for differentiating
under the integral sign}
\vskip.25in
Erik Talvila\\ [2mm]
\end{center}
{\bf 1. Introduction.}\quad
When we have an integral that depends on a parameter, say
$F(x)=\int_a^bf(x,y)\,dy$, it is often important to know
when $F$ is differentiable and when $F'(x)=\int_a^bf_1(x,y)\,dy$.
A sufficient condition for differentiating under the integral sign
is that $\int_a^bf_1(x,y)\,dy$ converges uniformly;  
see \cite[p.~260]{rogers}.
When we have absolute convergence, the condition $|f_1(x,y)|\leq g(y)$
with $\int_a^bg(y)\,dy<\infty$ suffices (Weierstrass M-test and
Lebesgue Dominated Convergence).  If we use the Henstock integral, then
it is not difficult to give necessary and sufficient conditions for
differentiating under the integral sign.  The conditions depend on
being able to integrate every derivative.

If $g\fn[a,b]\to\R$ is continuous on $[a,b]$ and  differentiable on $(a,b)$ 
it is not always the
case that $g'$ is Riemann or Lebesgue integrable over $[a,b]$.
However,
the Henstock integral integrates all derivatives and thus leads to the
most complete version of the Fundamental Theorem of Calculus.  
The Henstock
integral's definition in terms of Riemann sums is only slightly
more complicated than for the Riemann integral (simpler than the
improper Riemann integral), yet it includes the Riemann, improper Riemann,
and Lebesgue integrals as special cases.  Using the very strong version
of the Fundamental Theorem we can formulate necessary and sufficient conditions for
differentiating under the integral sign.

\bigskip
\noindent
{\bf 2. An introduction to the Henstock integral.}\quad
Here we lay out the facts about Henstock integration that we need.
There
are now quite a number of works that deal with this integral;  two good
ones to start with are
\cite{bartle1} and \cite{gordon}.

Let $f\fn[-\infty,\infty]\to(\infty,\infty)$.  A {\it gauge} is 
a mapping $\gamma$
from $[-\infty,\infty]$ to the open intervals in $[-\infty,\infty]$.  By 
{\it open
interval} we mean $(a,b)$, $[-\infty,b)$, $(a,\infty]$, or
$[-\infty,\infty]$ for
all $-\infty\leq a<b\leq \infty$ (the two-point compactification of the real 
line).
The defining property of the gauge is that for all $x\in[-\infty,\infty]$,
$\gamma(x)$ is an open interval containing $x$.  A {\it tagged partition} of
$[-\infty,\infty]$ is a finite set of pairs ${\cal P}=\{(z_i,I_i)\}_{i=1}^N$,
where
each
$I_i$ is a nondegenerate closed interval in $[-\infty,\infty]$ and $z_i\in I_i$.
The points $z_i \in[-\infty,\infty]$ are called {\it tags} and 
need not be distinct.
The intervals $\{I_i\}_{i=1}^N$ form a {\it partition}: For $i\not=j$, $I_i\cap I_j$ is empty
or a singleton and $\cup_{i=1}^NI_i=[-\infty,\infty]$.  We say $\cal P$ is
$\gamma$-{\it fine} if 
$I_i\subset\gamma(z_i)$ for all $1\leq i\leq N$.  Let $|I|$
denote the length of an interval with $|I|=0$ for an unbounded interval.  Then,
$f$ is {\it Henstock integrable}, and we write
$\int_{-\infty}^\infty f=A$, if there is a real number
$A$ such that for all $\epsilon>0$ there is a gauge function $\gamma$ 
such that if
${\cal P}=\{(z_i,I_i)\}_{i=1}^N$ is any $\gamma$-fine tagged partition of
$[-\infty,\infty]$ then
$$
\left|\lsum_{i=1}^Nf(z_i)\,|I_i|-A\right|<\epsilon.
$$
Note that $N$ is not fixed and the partitions can have any finite number of terms.
We can integrate over an interval $[a,b]\subset[-\infty,\infty]$ by multiplying
the integrand with the characteristic function $\chi_{[a,b]}$.

The more dramatically a function changes near a point $z$, the smaller $\gamma(z)$
becomes.  With the Riemann integral the intervals are made uniformly small.  Here
they are locally small.
A function is {\it Riemann integrable} on a finite
interval if and only if the gauge can be taken to assign intervals of
constant length.  It is not too surprising that the Henstock integral includes
the Riemann integral.  What is not so obvious is that the Lebesgue integral is
also included.  And, the Henstock integral can integrate functions that are
neither Riemann nor Lebesgue integrable.  An example is the function $f=g'$ where
$g(x)=x^2\sin(1/x^3)$ for $x\not=0$ and $g(0)=0$;
the origin is the only point of nonabsolute summability.  See
\cite[p.~148]{jeffrey} for a function that is Henstock integrable but whose
points of nonabsolute summability have positive measure.
A key feature of the Henstock integral is that it is nonabsolute:
an integrable function need not have an integrable absolute
value.

The convention $|I|=0$ for an unbounded interval performs essentially the same
truncation that is done with improper Riemann integrals and the Cauchy extension
of Lebesgue integrals.  A consequence of this is that there are
no improper Henstock integrals.
This fact is captured in the
following theorem,
which is proved for finite intervals in \cite{gordon}.
\begin{theorem}\label{improper}
Let $f$ be a real-valued function on $[a,b]\subseteq[-\infty,\infty]$.
Then $\int_a^bf$ exists and equals $A\in\R$ if and only if
$f$ is integrable on each subinterval $[a,x]\subset[a,b]$ and
$\lim_{x\to b^-}\int_a^xf$ exists and equals $A$.
\end{theorem}

Lebesgue integrals can be characterised by the fact that the indefinite integral
$F(x)=\int_a^xf$ is absolutely continuous.  A similar characterisation is
possible with the Henstock integral.  We need three definitions.  Let
$F\fn[a,b]\to\R$.  We say $F$ is {\it absolutely continuous} ($AC$) on a set
$E\subseteq[a,b]$ if for each $\epsilon>0$ there is some $\delta>0$ such that
$\sum_{i=1}^N|F(x_i)-F(y_i)|<\epsilon$ for all finite sets of disjoint open
intervals $\{(x_i,y_i)\}_{i=1}^N$ 
with endpoints in $E$ and $\sum_{i=1}^N(y_i-x_i)<\delta$.  We say that
$F$ is {\it absolutely continuous in the restricted sense} ($AC_*$) if instead
we have $\sum_{i=1}^N\sup_{x,y\in[x_i,y_i]}|F(x)-F(y)|<\epsilon$ under the
same conditions as with $AC$.  And, $F$ is said to be {\it generalised 
absolutely continuous 
in the restricted sense} ($ACG_*$) if $F$ is continuous and $E$ is the countable
union of sets on each of which $F$ is $AC_*$. 
Two useful properties are that among continuous functions, the $ACG_*$ 
functions
are properly contained in the class of functions that are 
differentiable almost everywhere
and they properly contain the class of functions that are differentiable 
nearly everywhere
(differentiable except perhaps on a countable set).  A function $f$ is
{\it Henstock integrable} if and only if there is an $ACG_*$  function $F$ with $F'=f$ almost
everywhere.  In this case $F(x) -F(a)=\int_a^xf$.  For an unbounded interval
such as $[0,\infty]$, continuity of $f$ at $\infty$ is obtained
by demanding that $\lim_{x\to\infty}F(x)$ exists.

We have given an example that shows that not
all derivatives are Lebesgue integrable.
However, all derivatives are
Henstock integrable.  This leads to a very strong version of the Fundamental
Theorem of Calculus.  A proof can be pieced together from results in \cite{gordon}.
\begin{theorem}[Fundamental
Theorem of Calculus]\mbox{}\\
{\bf I} Let $f\fn[a,b]\to\R$.  Then $\int_a^bf$ exists and $F(x)=\int_a^xf$
for all $x\in[a,b]$ if and only if $F$ is $ACG_*$ on $[a,b]$, $F(a)=0$, and
$F'=f$ almost everywhere on $(a,b)$.  If $\int_a^bf$ exists and $f$ is 
continuous at $x\in(a,b)$
then $\tfrac{d}{dx}\int_a^xf=f(x)$.\\
{\bf II} Let $F\fn[a,b]\to\R$.  Then $F$ is $ACG_*$ if and only if
$F'$ exists almost everywhere on $(a,b)$, $F'$ is Henstock integrable 
on $[a,b]$,
and $\int_a^xF'=F(x)-F(a)$ for all $x\in [a,b]$.
\end{theorem}
Here is a  useful sufficient condition for integrability of the 
derivative:
\begin{corollary}
Let $F\fn[a,b]\to\R$ be continuous on $[a,b]$ and differentiable nearly
everywhere on $(a,b)$.  Then $F'$ is Henstock integrable on $[a,b]$ and 
$\int_a^xF'=F(x)-F(a)$ for all $x\in [a,b]$.
\end{corollary}

The improvement over the Riemann and Lebesgue cases is that we need not
assume the integrability of $F'$.
Integration and differentiation are now inverse operations. To make this
explicit, let $A$ be the vector space of Henstock integrable functions on
$[a,b]\subseteq[-\infty,\infty]$, identified almost everywhere.  
Let $B$ be the vector space of
$ACG_*$ functions vanishing at $a$.  Let $\int$ be the integral operator defined
by $\int[f](x)=\int_a^xf$ for $f\in A$.  Let $D$ be the differential operator
defined by
$D[f](x)=f'(x)$ for $f\in B$.  The Fundamental Theorem then says that
$D\circ\int=I_A$ and $\int\circ D=I_B$.

\bigskip
\noindent
{\bf 3. Differentiation under the integral sign.}\quad 
\begin{theorem}
Let $f\fn[\alpha, \beta]\times[a,b]\to \R$.  Suppose that
$f(\cdot,y)$ is $ACG_*$ on $[\alpha,\beta]$ for almost all
$y\in(a,b)$.  Then
$F:=\int_a^bf(\cdot,y)\,dy$ is $ACG_*$ on $[\alpha,\beta]$
and $F'(x)=\int_a^bf_1(x,y)\,dy$ for almost all $x\in(\alpha,\beta)$ if and only if
\begin{equation}
\lint_{x=s}^t\lint_{y=a}^bf_1(x,y)\,dy\,dx=\lint_{y=a}^b\lint_{x=s}^tf_1(x,y)\,dx\,dy
\quad\mbox{for all } [s,t]\subseteq[\alpha,\beta].
\label{II.1}
\end{equation}
\end{theorem}

\bigskip
{\bf Proof:} Suppose $F$ is $ACG_*$ and $\frac{\partial}{\partial x}\int_{a}^bf(x,y)\,dy
=\int_{a}^bf_1(x,y)\,dy$.  Let $[s,t]\subseteq[\alpha,\beta]$.  By the second part of
the Fundamental Theorem, applied first to $F$ and then to $f(\cdot, y)$,
\begin{eqnarray}
\lint_{x=s}^t\lint_{y=a}^bf_1(x,y)\,dy\,dx & = & F(t)-F(s)\nonumber\\
 & = & \lint_{y=a}^b\left[f(t,y)-f(s,y)\right]\,dy\label{II.2}\\
 & = & \lint_{y=a}^b\lint_{x=s}^tf_1(x,y)\,dx\,dy.\nonumber
\end{eqnarray}

Now assume \eqref{II.1}.  Let $x\in(\alpha,\beta)$ and let $h\in\R$ be
such that
$x+h\in(\alpha,\beta)$.  Then, applying the second part of
the Fundamental Theorem to $f(\cdot, y)$ gives
\begin{eqnarray}
\lint_{x'=x}^{x+h}\lint_{y=a}^bf_1(x',y)\,dy\,dx' & = & 
\lint_{y=a}^b\lint_{x'=x}^{x+h}f_1(x',y)\,dx'\,dy\nonumber\\
 & = & \lint_{y=a}^b\left[
 f(x+h,y)-f(x,y)\right]\,dy\nonumber\\
 & = & \lint_{y=a}^bf(x+h,y)\,dy-\lint_{y=a}^bf(x,y)\,dy.\label{II.3}
 \end{eqnarray}
And,
\begin{eqnarray*}
F'(x) & = & \llim_{h\to 0}\tfrac1h\left[F(x+h)-F(x)\right]\\
 & = & \llim_{h\to 0}\tfrac1h\lint_{x'=x}^{x+h}\lint_{y=a}^bf_1(x',y)\,dy\,dx'\\
 & = & \int_a^bf_1(x,y)\,dy\quad\mbox{for almost all } x\in(\alpha,\beta).
\end{eqnarray*}
The last line comes from the first part of the Fundamental Theorem.
Repeating the argument in \eqref{II.2} shows that
$\int_\alpha^xF'=F(x)-F(\alpha)$ 
for all
$x\in[\alpha,\beta]$.  Hence, $F$ is $ACG_*$ on $[\alpha,\beta]$.\qed

The theorem holds for $-\infty\leq \alpha<\beta\leq\infty$ and 
$-\infty\leq a<b\leq\infty$.
Partial versions of the theorem are given in
\cite[p.~357]{hobson} for the Lebesgue integral
and in \cite[p.~63]{celidze} for the wide Denjoy integral, which includes the
Henstock integral.

From the proof of the theorem it is clear that only the linearity of
the integral over $y\in[a,b]$ (\eqref{II.2} and \eqref{II.3}) comes into play
with this variable.  Hence, we have the following generalisation.
\begin{corollary}\label{corII.1}
Let $S$ be some set and suppose $f\fn[\alpha, \beta]\times S\to \R$.  Let
$f(\cdot,y)$ be $ACG_*$ on $[\alpha,\beta]$ for all $y\in S$.
Let $T$ be the real-valued functions on $S$ and let 
${\cal L}$ be a linear functional defined on a subspace of $T$.
Define $F\fn [\alpha,\beta]\to\R$ by $F(x)={\cal L}[f(x,\cdot)]$.  Then $F$ 
is $ACG_*$ on $[\alpha,\beta]$ and
$F'(x)={\cal L}[f_1(x,\cdot)]$ for 
almost all $x\in(\alpha,\beta)$ if and only if
$$
\int_{s}^t{\cal L}[f_1(x,\cdot)]\,dx={\cal L}\int_{s}^tf_1(x,\cdot)\,dx
\quad\mbox{for all } [s,t]\subseteq[\alpha,\beta].
$$
\end{corollary}

If $f(x,y)=\int_\alpha^xg(x',y)\,dx'$ then $f(\cdot,y)$ is automatically
$ACG_*$.  This gives necessary and sufficient conditions for interchanging
iterated integrals.
\begin{corollary}\label{corII.2}
Let $g\fn[\alpha,\beta]\times[a,b]\to\R$.  Suppose that
$g(\cdot,y)$ is integrable over $[\alpha,\beta]$ for almost all
$y\in(a,b)$.  Define 
$G(x)=\int_a^b\int_\alpha^x
g(x',y)\,dx'\,dy$.  Then $G$ is $ACG_*$ on $[\alpha,\beta]$ and
$G'(x)=\int_a^bg(x,y)\,dy$ for almost all $x\in(\alpha,\beta)$ if and only if
\begin{equation}
\lint_{x=s}^t\lint_{y=a}^bg(x,y)\,dy\,dx=\lint_{y=a}^b\lint_{x=s}^tg(x,y)\,dx\,dy
\quad\mbox{for all } [s,t]\subseteq[\alpha,\beta].
\end{equation}
\end{corollary}

Combining Corollaries \ref{corII.1} and \ref{corII.2} 
gives necessary and sufficient conditions
for interchanging summation and integration.
\begin{corollary}\label{corII.3}
Let $g\fn[\alpha,\beta]\times\N\to\R$ and write $g_n(x)=g(x,n)$ for
$x\in[\alpha,\beta]$ and $n\in\N$.  Suppose that
$g_n$ is integrable over $[\alpha,\beta]$ for each $n\in\N$.  Define
$G(x)=\sum_{n=1}^\infty\int_\alpha^x
g_n(x')\,dx'$.  Then $G$ is $ACG_*$ on $[\alpha,\beta]$ and
$G'(x)=\sum_{n=1}^\infty g_n(x)$ for almost all $x\in(\alpha,\beta)$ if and only if
\begin{equation}
\lint_{x=s}^t\lsum_{n=1}^\infty g_n(x)\,dx=\lsum_{n=1}^\infty\,\lint_{x=s}^t
g_n(x)\,dx
\quad\mbox{for all } [s,t]\subseteq[\alpha,\beta].
\end{equation}
\end{corollary}

The Fundamental Theorem and its corollary yield conditions sufficient to allow
differentiation under the integral.
\begin{corollary}
Let $f\fn[\alpha, \beta]\times[a,b]\to \R$.\\
{\bf i)}\quad Suppose that
$f(\cdot,y)$ is continuous on $[\alpha,\beta]$ for almost all $y\in(a,b)$
and is differentiable
nearly everywhere in $(\alpha,\beta)$ for almost all $y\in(a,b)$.  
If \eqref{II.1}
holds then $F'(x)
=\int_{a}^bf_1(x,y)\,dy$ for almost all $x\in(\alpha,\beta)$.\\
{\bf ii)}\quad Suppose that
$f(\cdot,y)$ is $ACG_*$ on $[\alpha,\beta]$ for almost all
$y\in(a,b)$ and 
that $\int_a^bf_1(\cdot,y)\,dy$ is continuous on $[\alpha,\beta]$.
If \eqref{II.1} holds then $F'(x)=\int_{a}^bf_1(x,y)\,dy$ for all $x\in(\alpha,\beta)$.
\end{corollary}

Here is an example of what can go wrong when 
one differentiates under the integral
sign without justification.  In 1815 Cauchy obtained the convergent integrals
$$
\lint_{x=0}^{\infty}\left\{\!\!\begin{array}{c}
\sin(x^2)\\
\cos(x^2)
\end{array}
\!\!\right\}
\cos(sx)\,dx  = 
\frac{1}{2}\sqrt{\frac{\pi}{2}}\left[\cos\left(\frac{s^2}{4}\right)
\mp\sin\left(\frac{s^2}{4}\right)\right].
$$
He then differentiated under the integral sign with respect to
$s$ and obtained the two
divergent integrals 
\begin{eqnarray*}
\lint_{x=0}^{\infty}x\left\{\!\!\begin{array}{c}
\sin(x^2)\\
\cos(x^2)
\end{array}
\!\!\right\}
\sin(sx)\,dx & \mbox{`}=\mbox{'} & 
\frac{s}{4}\sqrt{\frac{\pi}{2}}\left[\sin\left(\frac{s^2}{4}\right)
\pm\cos\left(\frac{s^2}{4}\right)\right].
\end{eqnarray*}
These divergent integrals have been reproduced
ever since and still appear in standard tables today, listed as
converging.  This story
was told in \cite{talvila}.

{\it
University of Alberta,
Edmonton AB Canada T6G 2E2\\
etalvila@math.ualberta.ca}
\end{document}